\documentclass{article}

\usepackage[cp1251]{inputenc}
\usepackage[russian]{babel}
\usepackage{amsfonts, amsmath, amsthm}
\usepackage{graphicx}

\newtheorem{Th}{Теорема}
	\newtheorem{Lemm}{Лемма}
	\newtheorem{Example}{Пример}
	
	\newtheorem{St}{Предложение}[section]
		\theoremstyle{definition}
	
	\newtheorem{Def}{Определение}[section]
	
\newenvironment{Proof} 		
	{		
		\parskip=-3pt
		\parindent=12pt
		\par
		{\bf Доказательство.}
	}	
	{		
		\hfill
		$\square$
	}	
	
\newenvironment{Des}[1]		
	{
		{\bf Обозначения:}
		\par
		\vspace{5pt}
		\hbox{
			\hspace{10pt}
			\parbox{0.96\textwidth}%
			{#1}
		}
	}
	{
		\par
	}

\begin{document}

	\begin{center} 		
		\textbf{\Large Минимальное число запретов задающее периодическое слово} \\[15pt]
		{\large Лавров Петр}	\\[15pt]
		Московский Государственный Университет им. М. В. Ломоносова
	\end{center}
	
\section{Введение}

Пусть $A$ --- ассоциативная алгебра. Порядок
$a_1 < \cdots < a_s$ на ее образующих индуцирует порядок на множестве  мономов от
$a_i$ (сперва по длине, потом лексикографически). Множество мономов, не
являющихся линейной комбинацией меньших мономов, образует {\it нормальный
базис} алгебры $A$.  {\it Функция роста} $V_A(n)$ есть размерность
пространства, порожденного мономами степени не выше $n$.

Проблематика, связанная с канонической формой алгебр, заданных своим копредставлением (образующими и соотношениями) традиционно находится в центре внимания алгебраистов. Обширная библиография посвящена также проблемам роста и рациональности (см., например,  обзоры --- \cite{BBL,Latyshev,Semenov,Latyshev1}). Знаменитая работа М.~Л.~Громова была посвящена описанию групп полиномиального роста (\cite{Gromov}). 

А.~И.~Ширшов при решении проблемы равенства в алгебрах Ли с одним соотношением ввел понятие, которое в дальнейшем стало называться {\it базисом Грёбнера-Ширшова} или {\it базисом Грёбнера}. Рассмотрим старшие члены элементов идеала соотношений $I$ или {\it редуцируемые слова}. Заметим, что надслово редуцируемого слова --- редуцируемо. {\it Обструкцией} называется минимальное редуцируемое слово, т.е. не содержащее собственных редуцируемых подслов.
Если все обструкции содержатся во множестве старших членов некоторого базиса $I$, то этот базис называется {\it базисом Грёбнера-Ширшова} идеала $I$. Аналогичные понятия для полилинейных слов ввел В.~Н.~Латышев при получении оценок на функцию короста для $T$-идеалов (минимальность понимается также и в том, что слово не является изотонным образом редуцируемого). 
А.~И.~Ширшов ввел понятие {\it композиции} и предложил критерий того,
что $\{f_i\}$ является базисом Грёбнера. Это легло в основе знаменитой
Diamond-леммы Бергмана. Вычислению базисов Грёбнера для различных идеалов посвящено целое направление исследований.

В.~Н.~Латышев поставил вопрос об алгоритмической разрешимости проблемы распознавания делителей нуля при условии конечности базиса Грёбнера (решен отрицательно в работе \cite{Ivanov})  и для мономиального случая (алгоритмическая разрешимость показана в работах 
\cite{BBL,Belov}). А.~Я.~Белов поставил вопрос о росте минимального базиса Грёбнера или о {\it коросте} в алгебрах --- функции $B(n)$ равной числу обструкций длины не больше $n$. Если базис Грёбнера конечен, то корост --- константа. 
При изучении короста фундаментальное значение (в том числе для проблем Бернсайдовского типа) имеет вопрос о {\it кодлине} периодической последовательности периода $n$ или минимальном количестве запретов, задающих периодическую последовательность. Задачи о его нахождении были поставлены А.~Я.~Беловым в 1994 году. В 2007  году Г.~Р.~Челноковым было получено продвижение по данной задаче (\cite{Chelnokov}). 
Мы получили точную оценку для данного случая. Основной результат данной работы состоит в следующем: 

	\begin{Th}
Для произвольного алфавита для $n$-периодического слова заданного $k$ запретами выполняется оценка: $n\leq\varphi_k$, где $\varphi_k$ --- это $k$-е число Фибоначчи ($\varphi_0=\varphi_1=1$).
	\end{Th}

Приведенная оценка точна --- существует пример последовательности слов с соответствующим коростом:

\begin{Example}
$l_0 = "a"$, $l_1 = "b"$, $l_{i+1} = l_il_{i-1}$. ($l_2 = "ba"$, $l_3 = "bab" \cdots$) Тогда $(l_k)^\infty$ задается $k$ запретами.
\end{Example}

{\bf Замечания. \\
1.} Мы получаем экспоненциальный рост длины периода от размера системы запретов. \\
{\bf 2.} Конечная система запретов, задающая данное n-периодическое слово $\omega$, всегда существует --- например, все слова длины не больше $n+1$ не являющиеся подсловами данного слова. \\
{\bf 3.} Легко понять, что разные системы могут задавать одно и то же слово, и количество запретов не вполне строго определено. Однако оказывается, что все они совпадают с точностью до редукции и существует универсальная, так называемая \textit{приведенная} система запретов, задающая данное слово.

Полный текст доказательства см. \cite{MyBasic}.

В последствии  И.~И.~Богдановым и Г.~Р.~Челноковым другим методом был получен результат, аналогичный нашему для двусимвольного алфавита.\\

В основе нашего доказательства лежит работа с {\it графами и схемами Рози}.
Последовательности схем Рози исследовались в ряде работ. С их помощью
удается решить ряд алгоритмических проблем (см., например, \cite{Mitrofanov1}, \cite{Mitrofanov2}, \cite{Mitrofanov3}).
Отметим,что результаты работ \cite{Mitrofanov1} и \cite{Mitrofanov2} другим методом независимо получены
Ф.Дюрандом: \cite{Durand1}, \cite{Durand2}.
В терминах размеченных схем Рози удается получить критерий того, что
слово отвечает перекладыванию отрезков (\cite{BelovChernInt}).

{\bf Благодарности} Автор благодарен А.~Я.~Белову и А.~В.~Михалеву за внимание к работе и помощь, И.~А.~Рипсу за полезные обсуждения. Исследования были поддержаны грантом РФФИ №14-01-00548

\newpage

\section{Основные обозначения и определения}

Пусть $A = \{a_1,\cdots,a_m\}$ --- конечный $m$-буквенный алфавит.\\[-5pt]

\begin{Des}{
 I --- множество бесконечных в обе стороны слов в алфавите A заданных с
 точностью до сдвига, $I_n \subset I$ --- слова наименьшего периода n. 
\\
 $F$ --- множество конечных слов в алфавите A. $F_n \subset F$ --- слова длины n. $U_n = \bigcup\limits^{n}_{k=1} F_k$ --- слова длины не больше n.
\\
 $F^\omega = \{l \prec \omega, l \in F, \omega \in I\}$ --- конечные подслова $\omega$.\\
 $S^\omega = F \backslash F^\omega$ --- конечные неподслова $\omega$.\\}
	\end{Des}

Мы будем говорить о том, как можно задать, идентифицировать какое-то бесконечное периодическое слово из $I$ с помощью указания каких-то его обязательных неподслов, иначе говоря --- запретов. Для пояснения приведем пример: 

	\begin{Example} 
Периодическое слово $\cdots ababab\cdots = (ab)^\infty$ в алфавите $\{a,b\}$ можно задать двумя запретами --- $\{aa, bb\}$. Никаких других слов из $I$ в которых этих запрещенных подслов не встречается --- нет (после a всегда должна идти b, а после b --- a). Слово $(aabab)^\infty$ --- запретами ${aaa, bb, aabaa, babab}$, а слово $(a^{n-1}b)^\infty$ --- n запретами $\{a^n, b^2, ba^kb, 1\leq k < n-1\}$
	\end{Example}

Формально это можно описать так:

	\begin{Def} 		
		Будем называть любое подмножество $S\subseteq F$ {\it системой запретов,} а его элементы --- {\it запретами}. 	
		Будем говорить, что слово $\omega \in I$ удовлетворяет системе запретов S, если $\omega$ не содержит слов из S как подслова, или, что эквивалентно, $S \subseteq S^\omega$. 
		Будем говорить, что S задает $\omega$ ($S \propto \omega$), если оно единственное ей удовлетворяет.
	\end{Def}

Наименьший период характеризует сложность самого слова, а его кодлина или размер системы запретов --- сложность его описания. 




\begin{Def}
Будем называть систему запретов S задающую данное периодическое слово $\omega$ {\it приведенной}, если выполнены следующие условия:
1) S не содержит дубликатов
2) Запреты из S нельзя сократить: $\forall s \in S: l \prec s \& l \neq s \Rightarrow l \prec \omega (\Leftrightarrow l \notin S^\omega)$
\end{Def}

Такая система запретов существует и единственна. Любая система может быть заменена на эквивалентную (задающую то же слово) приведенную.

\section{Графы Рози}

Для оценки размера периода данного слова, размера его приведенной системы запретов и связи между ними мы будем строить некоторую последовательность графов, строение и размеры которых будут связаны с соответствующими параметрами нашего слова, так что через одно можно будет оценить другое, что мы и сделаем.

\begin{Def}
графом Рози с номером $k$ слова $\omega \in I$ называется граф $G_k^\omega$ такой что $V(G_k^\omega) = F^\omega \cap F_n$ и $E(G_k^\omega) = F^\omega \cap F_{n+1}$
\end{Def}

Эти графы естественным образом ориентированны

\begin{Example}
На рисунке 1 изображен нулевой граф любого нетривиального слова в двухбуквенном алфавите, а также 1-й и 2-й графы для слова $(aabab)^\infty$
\end{Example}

\includegraphics[height=25mm]{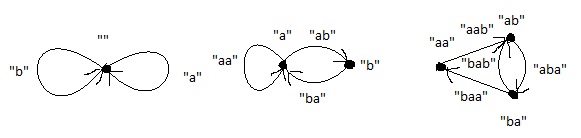} 


Рассмотрим как связаны два последовательных графа Рози одного и того же слова. Ясно, что $V(G_{k+1}^\omega) = E(G_k^\omega)$. Кроме того, каждому ребру графа $G_{k+1}^\omega$ соответствует путь длины два в $G_k^\omega$ --- слову из $k+2$ букв --- слово из его первых $k+1$ букв, и слово из его последних $k+1$ букв. Пересекаются они по средним $k$.

\begin{Example}
Ребру $aba$ на рисунке 1 будут соответствовать два полседовательных ребра --- $ab$ и $ba$.
\end{Example}

Обратное соответствие, однако, места не имеет, поскольку эволюция графов определяется системой запретов.

\begin{St}
Каждому пути длины 2 в $G_k^\omega$ однозначно соответствует либо ребро в  $G_{k+1}^\omega$, либо запрет из приведенной системы запретов задающей данное слово. 
\end{St}
\begin{Proof}
Если некоторому пути из подслов $l_1 = "a_1 \cdots a_{k+1}"$ и $l_2 = "a_2 \cdots a_{k+2}"$ не соответствует ребра, значит слова $"a_1 \cdots a_{k+2}"$ нет в $\omega$ как подслова. Пусть его нет и в системе запретов. Раз слова $l_1$ и $l_2$ встречаются в $\omega$, то что-то идет в нем до и после этих слов --- бесконечное слева $u l_1$ и справа $l_2 v$.\\[-10pt]

Значит слово $ua_1\cdots a_{k+2}v\in I$ не равное $\omega$, т.к. в $\omega$ не содержится $a_1\cdots a_{k+2}$ тоже удовлетворяет приведенной системе запретов $\tilde{S}^\omega$, т.к. если какой-то запрет в нем содержится, то он содержится в $u l_1$, $l_2 v$ или содержит $a_1\cdots a_{k+2}$, чего быть не может, т.к. система запретов приведенная --- противоречие.
\end{Proof}

Приведенная система запретов содержит в точности все такие слова, и значит она единственна. Будем обозначать приведенную систему запретов данного слова $\tilde{S^\omega}$


\section{Двухбуквенный алфавит}

Теперь перейдем к рассмотрению свойств графов Рози исключительно слов в 2-буквенном алфавите \{a, b\}

Поскольку букв только две, в графах Рози есть только четыре типа вершин в зависимости от их входящих и исходящих степеней --- (1, 1), (1, 2), (2, 1) и (2, 2). Будем называть вершины типа (1, 2) {\it исходящими развилками}, (2, 1) --- {\it входящими}, а (2, 2) --- {\it перекрестками}.

Эволюция происходит так: в центре каждого ребра ставим по вершине, соединем вершины в центрах последовательных ребер, и затем удаляем запрещенные --- см. рисунок.

\includegraphics[height=30mm]{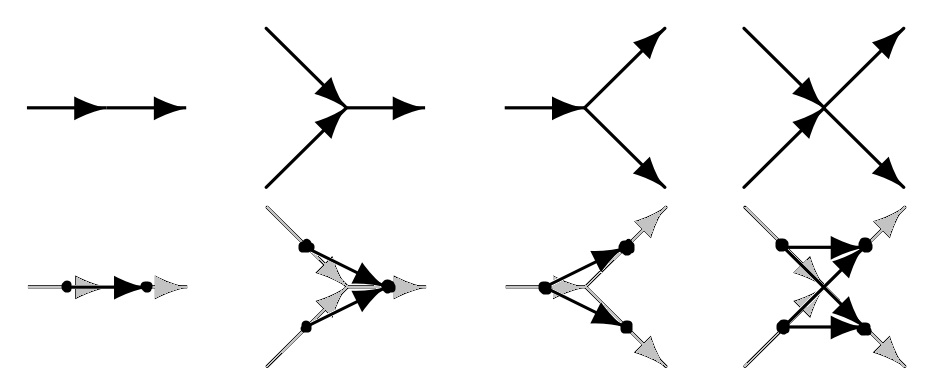}


Легко убедиться, что количество ребер графа $G^\omega_n$ n-периодического слова $\omega \in I_n$ в точности равно n (например, см. \cite{BBL}).


Теперь исследуем как связаны количества ребер, развилок и перекрестков в последовательных графах Рози. Введем обозначения:

 $i_k$ --- количество входящих развилок в $k$-м графе Рози 

 $c_k$ --- количество перекрестков 

	$inc_k$ --- общее количество вершин входящей степени 2

Сопоставим каждому пути длины 2 в $G_k$ его конечное ребро. Как мы выяснили, количество ребер в $G_{k+1}$ равно количеству таких путей за вычетом количества запретов соответствующего периода. Следовательно,
$$
|E(G_{k+1})| = |E(G_k)| + 2*c_k + i_k - |\tilde{S^\omega}\cap F_{n+2}|
$$
Здесь второе и третье слагаемое отвечают ребрам которые начинаются соответственно в перекрестках и входящих развилках, поскольку каждому из этих ребер соответствует два пути длины 2, а не один. Отсюда:
$$
n = |E(G_n)| = |E(G_0)| + \sum_{k=0}^{n-1} (|E(G_{k+1})| - |E(G_k)|) = |E(G_0)| + \sum_{k=0}^{n-1} (2*c_k + i_k - |\tilde{S^\omega}\cap F_{k+2}|)
$$



Посчитаем как меняется $inc_k$ (из равенства сумм входящих и исходящих степеней всех вершин очевидно, что их столько же сколько и вершин исходящей степени 2, также очевидно что $inc_k = i_k + c_k$). Для $G_{k+1}$ оно равно количеству ребер в $G_k$, начальная вершина которых имеет входящую степень 2, за исключением тех у которых эта степень уменьшится при запретах. Заметим, что для каждого удаленного ребра его конечная вершина обязана присутствовать в $\omega$ как подслово, поскольку она соответствовала конечному ребру некоторого двухреберного пути в $G_k$. Поэтому ее входящая степень после удаления минимум 1, а до удаления минимум 2, то есть ровно 2. Так что $inc_{k+1} = 2*c_k + i_k - |\tilde{S^\omega}\cap F_{n+2}|$.
То есть 
$$
inc_{k+1} - inc_k = c_k - |\tilde{S^\omega}\cap F_{n+2}|
$$

 Но в $G_0$ такая вершина одна, а в $G_n$ их 0. Заметим, что 
$$
inc_n = inc_0 + \sum_{k=0}^{n-1} (c_k - |\tilde{S^\omega}\cap F_{n+2}|)
$$
или  
$$
|\tilde{S^\omega}| = 1 + \sum_{k=0}^{n-1} c_k
$$

Следовательно, 
$$
|E(G_n)| = 2 + \sum_{k=0}^{n-1} (2*c_k + i_k - |\tilde{S^\omega}\cap F_{k+2}|) = 1 + \sum_{k=0}^{n-1} (c_k + i_k)
$$

Таким образом, мы получаем, что количество запретов равно общему количеству перекрестков во всех графах, а длинна конечного цикла примерно равна общему количеству развилок одного из типов во всех графах.


Теперь изучим эволюцию графов Рози на более глобальном уровне.

Отметим, что в последовательных графах легко установить соответствие между развилками, ведь каждой вершине входящей степени 2 соответствует одна или две таких вершины в следующем графе (если не было запретов) --- ребра, которые из нее выходили. Легко видеть из картинки (рис. 3), что если мы выберем некоторую пару развилок, длина пути между ними изменится одним из четырех способов. Т.е. если развилки "`смотрят"' друг на друга, то они приблизятся, а если в разные стороны --- удалятся.

\includegraphics[height=30mm]{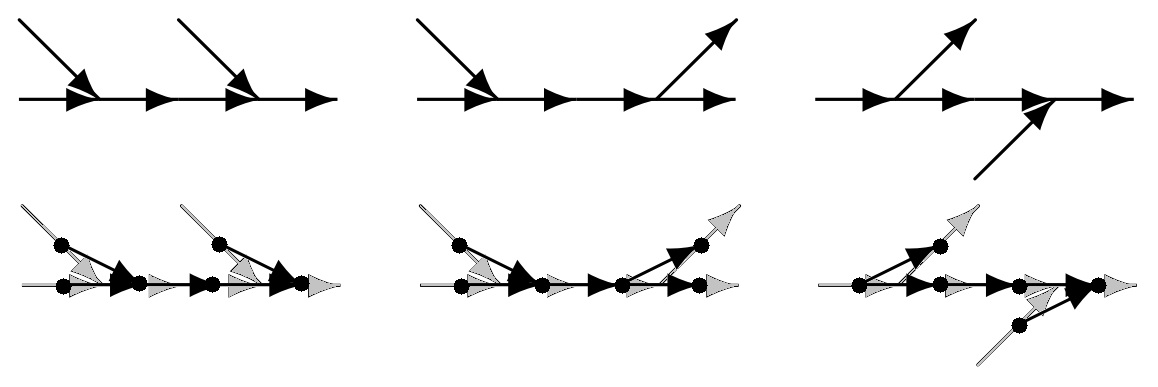}


Заметим, что в описанном случае происходит следующее --- развилки одновременно "`едут друг на друга"', "`сдвигаясь"' на $1/2$ ребра. 

Можно сделать по-другому --- остановить все развилки одного типа (допустим, входящие) и запустить развилки другого типа со скоростью 1. Тогда едущие развилки будут ``поглощать'' одно ребро и ``выдавать'' два --- сходно с тем как работает застежка молнии на куртке (рис. 4). 

\begin{picture}(40,40)
\put(0,0){\vector(1,1){15}}
\put(0,30){\vector(1,-1){15}}
\put(15,15){\vector(1,0){15}}
\put(30,15){\vector(1,0){15}}
\put(45,15){\vector(1,1){15}}
\put(45,15){\vector(1,-1){15}}
\put(60,30){\vector(1,0){15}}
\put(60,0){\vector(1,0){15}}

\put(90,0){\vector(1,1){15}}
\put(90,30){\vector(1,-1){15}}
\put(105,15){\vector(1,0){15}}
\put(120,15){\vector(1,1){15}}
\put(120,15){\vector(1,-1){15}}
\put(135,30){\vector(1,0){15}}
\put(135,0){\vector(1,0){15}}
\put(150,30){\vector(1,0){15}}
\put(150,0){\vector(1,0){15}}
\end{picture}

При этом, очевидно, граф будет меняться точно так же, и получится следующий граф Рози данного слова (более подробно этот момент разобран в полной статье, указанной в списке литературы). В конце пути развилки встречаются, и образуется перекресток. После этого, если не происходит запретов, исходящая развилка раздваивается и продолжает движение по обоим хвостам входящей. В случае запретов некоторые из образовавшихся развилок исчезают.

Получается, что образовавшемуся перекрестку естественным образом можно сопоставить длину пути который прошла исходящая развилка с момента своего появления до встречи с входящей. За время своей жизни эта развилка произвела именно такое количество ребер. Если бы каждая следующая развилка увеличивала суммарное количество раз не более чем в $\alpha$ раз, доказательство шло бы естественным прямолинейным путем. Однако это не так, и приходится делать поворот.

\section{Схема доказательства}


 Более общий способ получения экспоненциальной оценки --- создание некоторого упорядочивания развилок, альтернативного хронологическому, такого, что для него на каждом шаге (или нескольких шагах) сумма возрастает не более чем в соответствующее количество раз. 
Упорядочивание это строится с помощью следующего утверждения:



\begin{Lemm} \label{mainlemm}
Развилки можно двигать в произвольном порядке, не все параллельно, а последовательно, осуществляя те же запреты, что и в исходном слове. При этом количество запретов и конечное количество ребер не меняются.
\end{Lemm}

Отметим, что если промаркировать все развилки во всех исходных графах Рози, указать у них левую и правую сторону и запомнить какая маркировка устанавливается на образовавшихся развилках для каждой пары встретившихся развилок, а также какие ребра соответствующих развилок мы запрещаем и делать то же самое при не параллельном, а последовательном движении развилок, то в результате такой эволюции в конце получится также цикл, длина которого тоже равна периоду исходного слова. Это свойство в сущности является "`ассоциативностью"' для операций склеивания путей в графе (можно сказать, что в процессе эволюции некоторые пути между развилками склеиваются, чтобы в конце образовать цикл, и, очевидно, не важно в каком порядке их склеивать).

Будем "`схлапывать"' пары развилок между которыми ничего нет --- возьмем одну из этих развилок и подведем к другой "`расстегивая молнию"', таким образом формируя искомое упорядоченье --- соответсвенно порядку в котором производятся операции схлапывания (при которых образуется данный перекресток).

\begin{Lemm} \label{Estimation}
Если на некотором шаге для графа $G$ верна оценка $|E(G)| \leq \alpha^k * (\dfrac{2}{\alpha})^{d_G + 1}$, где $k$ --- номер текущего шага, а $d_G$ --- количество входящих развилок в $G$, либо $d_G = 1$ и пути в $G$ оцениваются числами Фибоначчи, то можно сделать несколько шагов так, чтобы указанная оценка снова оказалась верна.
\end{Lemm}


Можно показать, что размер конечного графа оценивается как двоичная экспонента от количества запретов - при каждой операции схлапывания размер графа будет увеличиваться на длинну схлапываемого пути. Путь не может быть длиннее всего графа --- а значит граф увеличиться не более чем в два раза на каждой операции.

Доказательство леммы \ref{Estimation} основанно на сходном соображении: если в графе $G$ одновременно $t > 2$ пар развилок, между которыми ничего нет --- можно "`схлопнуть"' их все. При этом размер графа также увеличится не более чем в два раза, а количество развилок уменьшится не более чем на $t$. Числовые вычисления показывают, что в таком случае оценка всегда сохраняется. Когда же $t \leq 2$ граф имеет достаточно простую структуру, поддающуюся детальному анализу. Оказывается что для этого случая выполнения оценки за несколько шагов также несложно добиться.

Многобуквенный случай доказывается с помощью того, что $r$-валентную развилку можно заменить на бинарное дерево с $r$ вершинами и для получившихся графов применить утверждение полученной для двухбуквенного случая оценки.

\end{document}